\documentclass[12pt]{article}
\usepackage[latin1]{inputenc}
\usepackage{amsmath, amsthm, amssymb}
\begin{document}

\bibliographystyle{plain}
\renewcommand{\Im}{\mathop{\rm Im }}
\renewcommand{\Re}{\mathop{\rm Re }}
\newcommand{\ra}{\mathop{\rightarrow }}
\newcommand{\supp}{\mathop{\rm supp}}
\newcommand{\sgn}{\mathop{\rm sgn }}
\newcommand{\card}{\mathop{\rm card }}
\newcommand{\KM}{\mbox{\rm KM}}
\newcommand{\diam}{\mathop{\rm diam}}
\newcommand{\diag}{\mathop{\rm diag}}
\newcommand{\tr}{\mathop{\rm tr}}
\newcommand{\Tr}{\mathop{\rm Tr}}
\newcommand{\dd}{ {\rm d} }
\newcommand{\id}{\mbox{\rm1\hspace{-.2ex}\rule{.1ex}{1.44ex}}
   \hspace{-.82ex}\rule[-.01ex]{1.07ex}{.1ex}\hspace{.2ex}}
\renewcommand{\P}{\mathop{\rm Prob}}
\newcommand{\V}{\mathop{\rm Var}}
\newcommand{\cps}{{\stackrel{{\rm p.s.}}{\longrightarrow}}}
\newcommand{\limm}{\mathop{\rm l.i.m.}}
\newcommand{\cloi}{{\stackrel{{\rm loi}}{\rightarrow}}}
\newcommand{\bra}{\langle\,}
\newcommand{\ket}{\,\rangle}
\newcommand{\obl}{/\!/}
\newcommand{\mapdown}[1]{\vbox{\vskip 4.25pt\hbox{\bigg\downarrow
  \rlap{$\vcenter{\hbox{$#1$}}$}}\vskip 1pt}}
\newcommand{\tab}{&\!\!\!}

\newcommand{\tabb}{&\!\!\!\!\!}
\renewcommand{\d}{\displaystyle}
\newcommand{\epreuve}{\hspace{\fill}$\bigtriangleup$}
\newcommand{\demo}{{\par\noindent{\em D\'emonstration~:~}}}
\newcommand{\solu}{{\par\noindent{\em Solution~:~}}}
\newcommand{\NB}{{\par\noindent{\bf Remarque~:~}}}
\newcommand{\const}{{\rm const}}
\newcommand{\cA}{{\cal A}}
\newcommand{\cB}{{\cal B}}
\newcommand{\cC}{{\cal C}}
\newcommand{\cD}{{\cal D}}
\newcommand{\cE}{{\cal E}}
\newcommand{\cF}{{\cal F}}
\newcommand{\cG}{{\cal G}}
\newcommand{\cH}{{\cal H}}
\newcommand{\cI}{{\cal I}}
\newcommand{\cJ}{{\cal J}}
\newcommand{\cK}{{\cal K}}
\newcommand{\cL}{{\cal L}}
\newcommand{\cM}{{\cal M}}
\newcommand{\cN}{{\cal N}}
\newcommand{\cO}{{\cal O}}
\newcommand{\cP}{{\cal P}}
\newcommand{\cQ}{{\cal Q}}
\newcommand{\cR}{{\cal R}}
\newcommand{\cS}{{\cal S}}
\newcommand{\cT}{{\cal T}}
\newcommand{\cU}{{\cal U}}
\newcommand{\cV}{{\cal V}}
\newcommand{\cW}{{\cal W}}
\newcommand{\cX}{{\cal X}}
\newcommand{\cY}{{\cal Y}}
\newcommand{\cZ}{{\cal Z}}
\newcommand{\bA}{{\bf A}}
\newcommand{\bB}{{\bf B}}
\newcommand{\bC}{{\bf C}}
\newcommand{\bD}{{\bf D}}
\newcommand{\bE}{{\bf E}}
\newcommand{\bF}{{\bf F}}
\newcommand{\bG}{{\bf G}}
\newcommand{\bH}{{\bf H}}
\newcommand{\bI}{{\bf I}}
\newcommand{\bJ}{{\bf J}}
\newcommand{\bK}{{\bf K}}
\newcommand{\bL}{{\bf L}}
\newcommand{\bM}{{\bf M}}
\newcommand{\bN}{{\bf N}}
\newcommand{\bP}{{\bf P}}
\newcommand{\bQ}{{\bf Q}}
\newcommand{\bR}{{\bf R}}
\newcommand{\bS}{{\bf S}}
\newcommand{\bT}{{\bf T}}
\newcommand{\bU}{{\bf U}}
\newcommand{\bV}{{\bf V}}
\newcommand{\bW}{{\bf W}}
\newcommand{\bX}{{\bf X}}
\newcommand{\bY}{{\bf Y}}
\newcommand{\bZ}{{\bf Z}}
\newcommand{\bu}{{\bf u}}
\newcommand{\bv}{{\bf v}}
\newfont{\msbm}{msbm10 scaled\magstep1}
\newfont{\msbms}{msbm7 scaled\magstep1} 
\newcommand{\bbA}{\mbox{$\mbox{\msbm A}$}}
\newcommand{\bbB}{\mbox{$\mbox{\msbm B}$}}
\newcommand{\bbC}{\mbox{$\mbox{\msbm C}$}}
\newcommand{\bbD}{\mbox{$\mbox{\msbm D}$}}
\newcommand{\bbE}{\mbox{$\mbox{\msbm E}$}}
\newcommand{\bbF}{\mbox{$\mbox{\msbm F}$}}
\newcommand{\bbG}{\mbox{$\mbox{\msbm G}$}}
\newcommand{\bbH}{\mbox{$\mbox{\msbm H}$}}
\newcommand{\bbI}{\mbox{$\mbox{\msbm I}$}}
\newcommand{\bbJ}{\mbox{$\mbox{\msbm J}$}}
\newcommand{\bbK}{\mbox{$\mbox{\msbm K}$}}
\newcommand{\bbL}{\mbox{$\mbox{\msbm L}$}}
\newcommand{\bbM}{\mbox{$\mbox{\msbm M}$}}
\newcommand{\bbN}{\mbox{$\mbox{\msbm N}$}}
\newcommand{\bbO}{\mbox{$\mbox{\msbm O}$}}
\newcommand{\bbP}{\mbox{$\mbox{\msbm P}$}}
\newcommand{\bbQ}{\mbox{$\mbox{\msbm Q}$}}
\newcommand{\bbR}{\mbox{$\mbox{\msbm R}$}}
\newcommand{\bbS}{\mbox{$\mbox{\msbm S}$}}
\newcommand{\bbT}{\mbox{$\mbox{\msbm T}$}}
\newcommand{\bbU}{\mbox{$\mbox{\msbm U}$}}
\newcommand{\bbV}{\mbox{$\mbox{\msbm V}$}}
\newcommand{\bbW}{\mbox{$\mbox{\msbm W}$}}
\newcommand{\bbX}{\mbox{$\mbox{\msbm X}$}}
\newcommand{\bbY}{\mbox{$\mbox{\msbm Y}$}}
\newcommand{\bbZ}{\mbox{$\mbox{\msbm Z}$}}
\newcommand{\bbsA}{\mbox{$\mbox{\msbms A}$}}
\newcommand{\bbsB}{\mbox{$\mbox{\msbms B}$}}
\newcommand{\bbsC}{\mbox{$\mbox{\msbms C}$}}
\newcommand{\bbsD}{\mbox{$\mbox{\msbms D}$}}
\newcommand{\bbsE}{\mbox{$\mbox{\msbms E}$}}
\newcommand{\bbsF}{\mbox{$\mbox{\msbms F}$}}
\newcommand{\bbsG}{\mbox{$\mbox{\msbms G}$}}
\newcommand{\bbsH}{\mbox{$\mbox{\msbms H}$}}
\newcommand{\bbsI}{\mbox{$\mbox{\msbms I}$}}
\newcommand{\bbsJ}{\mbox{$\mbox{\msbms J}$}}
\newcommand{\bbsK}{\mbox{$\mbox{\msbms K}$}}
\newcommand{\bbsL}{\mbox{$\mbox{\msbms L}$}}
\newcommand{\bbsM}{\mbox{$\mbox{\msbms M}$}}
\newcommand{\bbsN}{\mbox{$\mbox{\msbms N}$}}
\newcommand{\bbsO}{\mbox{$\mbox{\msbms O}$}}
\newcommand{\bbsP}{\mbox{$\mbox{\msbms P}$}}
\newcommand{\bbsQ}{\mbox{$\mbox{\msbms Q}$}}
\newcommand{\bbsR}{\mbox{$\mbox{\msbms R}$}}
\newcommand{\bbsS}{\mbox{$\mbox{\msbms S}$}}
\newcommand{\bbsT}{\mbox{$\mbox{\msbms T}$}}
\newcommand{\bbsU}{\mbox{$\mbox{\msbms U}$}}
\newcommand{\bbsV}{\mbox{$\mbox{\msbms V}$}}
\newcommand{\bbsW}{\mbox{$\mbox{\msbms W}$}}
\newcommand{\bbsX}{\mbox{$\mbox{\msbms X}$}}
\newcommand{\bbsY}{\mbox{$\mbox{\msbms Y}$}}
\newcommand{\bbsZ}{\mbox{$\mbox{\msbms Z}$}}

%
\def\eurtoday{\number\day \space\ifcase\month\or
 January\or February\or March\or April\or May\or June\or
 July\or August\or September\or October\or November\or December\fi
 \space\number\year}
%
%
\def\aujourdhui{\number\day \space\ifcase\month\or
 janvier\or f{\'e}vrier\or mars\or avril\or mai\or juin\or
 juillet\or ao{\^u}t\or septembre\or octobre\or novembre\or d{\'e}cembre\fi
 \space\number\year}
\textheight=21cm
\textwidth=16cm 
\voffset=-1cm
\newtheorem{theo}{Theorem}[section]
\newtheorem{pr}{Proposition}[section]
\newtheorem{cor}{Corollary}[section]
\newtheorem{lem}{Lemma}[section]
\newtheorem{defn}{Definition}[section]
\hoffset=-1,5cm
\parskip=4mm
\title{Discrete approximation of a stable self-similar stationary increments process}
\author{C. Dombry  and N. Guillotin-Plantard \protect\hspace{1cm}}
\date{}
\maketitle
~
~\\
Running title: Discrete approximation of SSSSI process.\\
~\\ 
Authors' addresses:\\
C. Dombry (corresponding author), Université de Poitiers, Laboratoire de Mathématiques et Applications, Téléport 2 - BP 30179,
Boulevard Marie et Pierre Curie, 86962 Futuroscope Chasseneuil Cedex, France, e-mail: clement.dombry@math.univ-poitiers.fr\\
N. Guillotin-Plantard, Universit\'e de Lyon, Institut Camille Jordan, 43 bld du 11 novembre 1918, 69622 Villeurbanne, France, e-mail: nadine.guillotin@univ-lyon1.fr \\

\begin{abstract}
The aim of this paper is to present a result of discrete approximation of some class of stable self-similar stationary increments processes. The properties of such processes were intensively investigated, but little is known on the context in which such processes can arise. To our knowledge, discretisation and convergence theorems are available only in the case of stable Lévy motions and fractional Brownian motions. This paper yields new results in this direction. \\
Our main result is the convergence of the random rewards schema, which was firstly introduced by Cohen and Samorodnitsky, and that we consider in a more general setting.
Strong relationships with Kesten and Spitzer's random walk in random sceneries are evidenced. Finally, we study some path properties of the limit process.  
\end{abstract}

\section{Introduction} 
\subsection{Motivations}
In the past few decades, many efforts were done to achieve a better comprehension of stable self-similar stationary increments processes, and to achieve a classification of such processes . The Gaussian case was studied by Mandelbrot and Van Ness \cite{MVN68} and gives rise to an easy classification: to any index of self-similarity $H\in(0,1)$ corresponds a unique (up to a multiplicative constant) Gaussian H-self-similar process, called {\it fractional Brownian motion}. Basic properties of the fractional Brownian motion were investigated by Mandelbrot and Van Ness. In the $\alpha$-stable case, $0<\alpha<2$, the family of self-similar stationary increments processes is much larger. The feasible range of pairs $(\alpha,H)$ (where $\alpha$ is the index of stability and $H$ the exponent of self-similarity) is 
$$\left\{ \begin{array}{ll} 0<H\leq \frac{1}{\alpha}& {\rm \ if\ } 0<\alpha \leq 1 \\ 0<H<1& {\rm \ if\ } 1<\alpha <2 \end{array}\right.$$
and, apart from the case $0<\alpha<1$ and $H=1/\alpha$, a feasible pair $(\alpha,H)$ does not determine the law of the process.
Two most well known families of stable self-similar stationary increments process are obtained by taking two of the many possible integral representations of the fractional Brownian motion and modifying them appropriately: these are the linear fractional stable motion and the real harmonizable fractional stable motion.\\
In many of these works, the construction of such processes is based on the theoretical tool of stochastic or stable integration. This strategy allows to construct many processes with interesting properties. However, apart from the case of Lévy motions and fractional Brownian motions, we do not know any situation where these processes naturally arise. For example, Donsker theorem yields an approximation of Brownian motion by discrete processes that arise naturally in many cases of modelisation, and hence explains why Brownian motion is so important. In the context of communication networks, the best known  result of this type is perhaps the paper by Mikosh et al. \cite{MRRS02}, where the limiting model turns out to be either fractional Brownian motion or Lévy stable motion (depending on the relationship between the number of users and the time scale). \\
In \cite{CS05}, the authors present a discrete schema, called the {\it random rewards schema}, and show its convergence to a stable self-similar statinonary increments process called FBM-1/2-local time fractional symmetric stable motion. This provides a rigorous framework, in which a non standard self-similar stationary increments process appears naturally as the limit of a discrete schema. 
Roughly speaking, the schema is the following: a simple network is modelised by the integer line; we suppose that to each node of the network is associated a random reward; a user is moving randomly on the network, and is earning the rewards associated to the nodes he visits; the schema deals with the total reward earned when many such users are considered.\\
In this paper, we consider the random rewards schema in a more general setting. We suppose that the random reward is in the domain of attraction of a $\beta$-stable distribution and that the random walk performed by the user is in the domain of attraction of a $\alpha$-stable Lévy motion, with $1<\alpha \leq 2$ (in \cite{CS05}, only the case $\alpha=2$ was considered.) The convergence of the random reward schema is proven and some properties of the limit process are investigated. This yields a new class of $\beta$-stable  $H$-self-similar stationary increments processes, with $H=1-\alpha^{-1}+(\alpha\beta)^{-1}$. With this procedure, the range of feasible pairs is :  
$$\left\{ \begin{array}{ll} \frac{\beta+1}{2\beta}\leq H< \frac{1}{\beta} & {\rm \ if\ } 0<\beta < 1 \\
H=1 &{\rm \ if\ } \beta=1 \\
 \frac{1}{\beta}<H\leq \frac{\beta+1}{2\beta} &{\rm \ if\ } 1<\beta <2. \end{array}\right.$$
The next section is devoted to the presentation of the model and of the results we obtained. In the second section we prove the convergence of the generalised random reward schema. We begin the proof with some estimations on quantities associated with random walks (local times, number of self intersections, range), and then focus on the convergence of the finite-dimensional distributions and on the tightness of the sequence.
In the last section, some properties of the limit process are studied.

\subsection{Model and results}
To place our results in context, we first describe a result of Kesten and Spitzer \cite{KS79} which shows that a stable process in random scenery can be realized as the limit in distribution of a random walk in random scenery.\\

Let $\xi=(\xi_x)_{x\in\bbsZ}$ denote a sequence of independent, identically distributed, real-valued random variables. 
The sequence $\xi$ is called a {\it random scenery}. Suppose that it belongs to the normal domain of attraction  of a strictly stable distribution $Z_{\beta}$ of index $\beta\in (0,2]$. This means that the following weak convergence holds:
\begin{equation}\label{eq1.01}
n^{-\frac{1}{\beta}}\sum_{x=0}^n \xi_x \mathop{\Longrightarrow}_{n\rightarrow\infty}^{\cL} Z_\beta. 
\end{equation}
We suppose that the characterisic function $\bar\lambda$ of the distribution $Z_\beta$ writes , for all $u\in\bbR$
\begin{equation}\label{eq1.02}
\bar\lambda(u)=\bbE\exp(iuZ_\beta)=\exp\left(-\sigma^\beta|u|^\beta(1-i\nu\tan(\frac{\pi\beta}{2})\sgn(u))\right), 
\end{equation}
for some constants $\sigma>0$ and $\nu\in [-1,1]$. Note that in the case $\beta\neq 1$, this is the most general form of a strictly stable distribution. In the case $\beta=1$, this is the general form of a strictly stable distribution satisfying the symmetry condition 
\begin{equation}\label{eq1.03} 
\exists M\geq 0 \ \ \forall L>0 \ \ \bbE(|Z_\beta |1_{\{|Z_\beta|<L\}}) \leq M .
\end{equation}

Let $S=(S_k)_{k\in\bbsN}$ be a {\it random walk} on $\bbZ$ independent of the random scenery $\xi$. We suppose that $$\left\{\begin{array}{l}S_0=0, \\ 
S_n=\sum_{k=1}^n X_k\ ,\ n\geq 1,\end{array}\right.$$ 
where the $X_i$ are i.i.d. integer-valued random variables belonging to the normal domain of attraction of a strictly stable distribution $Z_{\alpha}$ of index $\alpha\in (1,2]$. This means that the following weak convergence holds:
\begin{equation}\label{eq1.04}
n^{-\frac{1}{\alpha}} S_n \mathop{\Longrightarrow}_{n\rightarrow\infty}^{\cL} Z_\alpha, 
\end{equation}
where the distribution $Z_\alpha$ has characteristic function
\begin{equation}\label{eq1.05}
\bbE\exp(iuZ_\alpha)=\exp\left(-\chi^\alpha|u|^\alpha(1-i\mu\tan(\frac{\pi\alpha}{2})\sgn(u))\right), 
\end{equation}
for all $u\in\bbR$ and for some constants $\chi>0$ and $\mu\in [-1,1]$.

We define the {\it random walk in random scenery} as the process $(Z_n)_{n\geq 0}$ given by
\begin{equation}\label{eq1.06}
Z_n=\sum_{k=0}^n \xi_{S_k}.
\end{equation}
Stated simply, a random walk in random scenery is a cumulative sum process whose summands are drawn from the scenery; the order in which the summands are drawn is determined by the path of the random walk.
We extend this definition to non-integer time $s\geq 0$ by the linear interpolation
\begin{equation}\label{eq1.07}
Z_s=Z_{[s]}+(s-[s])(Z_{[s]+1}-Z_{[s]}).
\end{equation}

We now describe the limit theorem for the random walk in random scenery established by Kesten and Spitzer \cite{KS79}.\\ 
Cumulative sums of the scenery converge in $\cD(\bbR)$, the space of càd-làg functions:
$$\left(n^{-\frac{1}{\beta}} \sum_{k=0}^{[nx]}\xi_k\right)_{x\in \bbsR} \mathop{\Longrightarrow}_{n\rightarrow\infty}^{\cL} \left(W(x)\right)_{x\in \bbsR},$$ 
where $W$ is a bilater $\beta$-stable Lévy process such that $W(0)=0$, and $W(1)$ and $W(-1)$ are distributed according to  $Z_\beta$.\\
The random walk converge in $\cD([0,\infty))$: 
$$\left(n^{-\frac{1}{\alpha}} S_{[nt]}\right)_{t\geq 0} \mathop{\Longrightarrow}_{n\rightarrow\infty}^{\cL} \left(Y(t)\right)_{t\geq 0},$$
where the limit process $Y$ is an $\alpha$-stable Lévy process such that $Y(0)=0$ and $Y(1)$ is distributed according to $Z_\alpha$.\\ 
To describe the limit process known as {\it stable process in random scenery} we suppose that $Y$ and $W$ are two independent processes defined on the same probability space and distributed as above. Let $L_t(x)$ the jointly continuous version of the local time of the process $Y$.
Kesten and Spitzer prove the following weak convergence in the space of continuous function $\cC([0,\infty))$  
\begin{equation}\label{eq1.08}
\left(n^{-\delta} Z_{nt}\right)_{t\geq 0} \mathop{\Longrightarrow}_{n\rightarrow\infty}^{\cL}  \left(\Delta(t)\right)_{t\geq 0}
\end{equation} 
where $\delta=1-\alpha^{-1}+(\alpha \beta)^{-1}$ and $\Delta$ is the process defined by
$$\Delta(t)=\int_{-\infty}^{+\infty} L_t(x){\rm d}W(x).$$
The limit process $\Delta$ is known to be a  continuous $\delta$-self-similar stationary increments process. It can be seen as a mixture of $\beta$-stable processes, but it is not a stable process.
 
To get a stable process, we consider sums of independent stable processes in random scenery. 
Let $\Delta^{(i)}, i\geq 1$ be independent copies of the process $\Delta$. 
\begin{theo}\label{theo1} 
The following weak convergence holds in $\cC([0,\infty))$:
\begin{equation}\label{eq1.09}
\left(n^{-\frac{1}{\beta}} \sum_{i=1}^n \Delta^{(i)}(t)\right)_{t\geq 0} \mathop{\Longrightarrow}_{n\rightarrow\infty}^{\cL}  \left(\Gamma(t)\right)_{t\geq 0},
\end{equation}
\end{theo}
\noindent
{\bf Remark: }\ In this theorem, $\Gamma$ is defined as the weak limit of a sequence of processes. As we will see in Theorem \ref{theo3} below, $\Gamma$ is a continuous $\beta$-stable $\delta$-self-similar stationary increments process. It is worth noting that in the case $\beta=2$ (i.e. if the scenery has finite variance), $\Gamma$ must hence be a fractional Brownian motion of index $\delta=1-\frac{1}{2\alpha}$.

Replacing the stable process in random scenery by a random walk in random scenery, we obtain the random rewards schema which yields a discrete approximation of the process $\Gamma$.
Let $\xi^{(i)}=(\xi^{(i)}_x)_{x\in \bbsZ},\ i\geq 1$ be independent copies of $\xi$. Let $S^{(i)}=(S^{(i)}_n)_{n\in\bbsN}$ be independent copies of $S$ and also independent of the $\xi^{(i)}, i\geq 1$.
Denote by $D_n^{(i)}$ the $i$-th random walk in random scenery defined by 
\begin{equation}\label{eq1.10}
D_n^{(i)}(t)=n^{-\delta} Z_{nt}^{(i)} 
\end{equation}
where the definition of $Z_n^{(i)}$ is given by equations (\ref{eq1.06}) and (\ref{eq1.07}) with $\xi$ and $S$ replaced by the $i$-th random scenery $\xi^{(i)}$ and the $i$-th random walk $S^{(i)}$ respectively. 
\begin{theo}\label{theo2}
Let $c_n$ be a sequence of integers such that $\lim c_n=+\infty$. Then, the following weak convergence holds in $\cC([0,\infty))$:
\begin{equation}\label{eq1.11}
\left(c_n^{-\frac{1}{\beta}} \sum_{i=1}^{c_n} D_n^{(i)}(t)\right)_{t\geq 0} \mathop{\Longrightarrow}_{n\rightarrow\infty}^{\cL}  \left(\Gamma(t)\right)_{t\geq 0}.
\end{equation}
The limit process $\Gamma$ is the same as in Theorem \ref{theo1}.
\end{theo}
 
We then focus on the limit process $\Gamma$. Its finite dimensional distribution are $\beta$-stable multivariate distribution whose characteristic function are explicited in the following Theorem. Since the finite dimensional distributions are not very convenient to describe trajectorial properties of the process $\Gamma$, we also introduce an integral representation of $\Gamma$. \\
Let $(\Omega',\cF',\bbP')$ be a probability space supporting an $\alpha$-stable Lévy motion $Y(t)$ and let $L_t(x)$ be its joinly continuous local time process. Let $M$ be a $\beta$-stable random measure on the space $\Omega'\times\bbR$ with control measure $\bbP'\times\lambda$ (where $\lambda$ is the Lebesgue measure on $\bbR$) and skewness intensity identically equal to $\nu$. The random measure itself lives on some other probability space $(\Omega,\cF,\bbP)$. For the definition and properties of stable random measures and integrals with respect with these measures, we refer to the monography of Samorodnitsky and Taqqu \cite{ST94}.
\begin{theo}\label{theo3}
The process $\Gamma$ is a $\beta$-stable $\delta$-self-similar stationary increments process
whose finite dimensional distributions satisfy:\\ 
for all $(\theta_1,\cdots,\theta_k)\in \bbR^k$ and $(t_1,\cdots,t_k)\in[0,+\infty)^{k}$, 
\begin{eqnarray}
& &\bbE\left[\exp\left(i\sum_{j=1}^k \theta_j \Gamma(t_j)\right) \right]\label{eq1.12}\\
&=&\exp\left[-\sigma^\beta \bbE \int_{-\infty}^{+\infty} |\sum_{j=1}^k \theta_jL_{t_j}(x)|^\beta \left(1-i\nu \tan(\frac{\pi\beta}{2})\sgn(\sum_{j=1}^k \theta_jL_{t_j}(x))\right){\rm d}x  \right].\nonumber
\end{eqnarray}
Furthermore, it has the following integral representation:
\begin{equation}\label{eq1.13}
\Gamma(t)\mathop{=}^{\cL}\int_{\Omega'\times\bbsR} \sigma L_t(x)(\omega')M({\rm d}\omega',{\rm d}x)\ ,\ t\geq 0
\end{equation}
where $\overset{\cL}{=}$ stands for the equality of the finite dimensional distributions of the two processes.
\end{theo}

This integral representation is a powerful tool to derive the path properties of the process $\Gamma$. It can be used to explicit the probability tails of the suprema of the process. In the following, we only consider the case $0<\beta<2$. As we mentioned above, in the case $\beta=2$, the process $\Gamma$ is a fractional Brownian motion and its properties were widely studied. 
\begin{theo}\label{theo4}
Suppose that $0<\beta<2$. Then, for any $T>0$,
$$\lim_{u\rightarrow \infty} u^{\beta}\bbP\left(\sup_{t\in [0,T]} \Gamma(t)\geq u \right)
 =C_\beta\frac{1+\nu}{2} \bbE\left(\int_{\bbsR}L_T(x)^\beta {\rm d}x\right)$$
and
$$\lim_{u\rightarrow \infty} u^{\beta}\bbP\left(\sup_{t\in [0,T]} |\Gamma(t)|\geq u \right)
=C_\beta\bbE\left(\int_{\bbsR}L_T(x)^\beta {\rm d}x\right),$$
with 
$$C_\beta=\left(\int_0^\infty x^{-\beta}\sin x{\rm d}x\right)^{-1}.$$
\end{theo}

The last result of this paper is devoted to the Hölder regularity of sample path of the process $\Gamma$:
\begin{theo}\label{theo5}
Suppose that $0<\beta<2$. If $\beta \geq 1$, we suppose furthermore that $\nu=0$. Then, the process $\left(\Gamma(t)\right)_{t\geq 0}$ satisfies almost surely 
$$\sup_{0\leq s<t\leq 1/2}\frac{|\Gamma(t)-\Gamma(s)|}{(t-s)^{1-\frac{1}{\alpha}}|\log(t-s)|^{\frac{1}{\alpha}+\varepsilon}}<\infty,$$
with 
$$\varepsilon=\left\{ \begin{array}{ll} 0 &{\rm \ if\ } 1\leq \beta< 2 \\
\frac{1}{2} &{\rm \ if\ } 0< \beta< 1 \end{array}\right..$$
\end{theo}

\section{Sums of stable processes in random scenery \\ Proof of Theorem \ref{theo1}}
For $n\geq 1$, let $\Gamma_n$ the continuous process defined by 
$$\Gamma_n(t)=n^{-\frac{1}{\beta}}\sum_{i=1}^n \Delta^{(i)}(t)\ ,\ t\geq 0.$$
Theorem \ref{theo1} claims that the sequence $\Gamma_n$ converges weakly in $\cC([0,\infty))$. 
We prove this fact by proving the convergence of the finite dimensional distributions and the tightness of the sequence.
Theorem \ref{theo1} is thus a consequence of Propositions \ref{pr1.1} and \ref{pr1.2} below. 

\begin{pr}\label{pr1.1}
The finite dimensional distributions of $(\Gamma_n(t))_{t\geq 0}$ converge weakly as $n\rightarrow\infty$. 
\end{pr} 
{\it Proof :} \ Let $(\theta_1,\cdots,\theta_k)\in \bbR^k$ and $(t_1,\cdots,t_k)\in [0,+\infty)^k$. 
We prove the convergence of the characteristic functions
$$\bbE\left[\exp\left( i\sum_{j=1}^k \theta_j\Gamma_n(t_j)\right)\right]=\bbE\left[\exp\left( in^{-\frac{1}{\beta}}\sum_{i=1}^n\sum_{j=1}^k \theta_j \Delta^{(i)}(t_j)\right)\right].$$
Since the $\Delta^{(i)}$ are independent copies of $\Delta$,
\begin{equation}\label{eq2.1}
\bbE\left[\exp\left( i\sum_{j=1}^k \theta_j\Gamma_n(t_j)\right)\right]=\bbE\left[\exp\left( in^{-\frac{1}{\beta}}\sum_{j=1}^k \theta_j \Delta(t_j) \right)\right]^n.
\end{equation}
Using Lemma 5 in \cite{KS79} giving the characteristic function of the finite dimensional distributions of $\Delta$, we compute
\begin{equation}\label{eq2.2}
\bbE\left[\exp\left( in^{-\frac{1}{\beta}}\sum_{j=1}^k \theta_j \Delta(t_j) \right)\right]=\bbE\left[\exp(-n^{-1}X)\right]
\end{equation}
with
\begin{equation}\label{eq2.3}
X=\sigma^{\beta} \int_{\bbsR} |\sum_{j=1}^k \theta_jL_{t_j}(x)|^\beta\left(1- i\nu \tan(\frac{\pi\beta}{2})\sgn(\sum_{j=1}^k \theta_jL_{t_j}(x))\right) {\rm d}x .
\end{equation}
We will prove that equation (\ref{eq2.2}) implies the following asymptotic as $n\rightarrow\infty$ 
\begin{equation}\label{eq2.4}
\bbE\left[\exp\left( in^{-\frac{1}{\beta}}\sum_{j=1}^k \theta_j \Delta(t_j) \right)\right]=1-n^{-1}\bbE(X)+o(n^{-1}).
\end{equation}
The fact that $X$ is integrable is proven in Lemma \ref{lem2.1} below. We now prove equation (\ref{eq2.4}). To this aim, we
compute the limit as $n$ goes to infinity of
$$n\left(\bbE\left[\exp\left( in^{-\frac{1}{\beta}}\sum_{j=1}^k \theta_j \Delta(t_j) \right)\right]-1\right)=\bbE(f_n(X)),$$
where $f_n$ is defined on $\bbC$ by $f_n(x)=n(\exp(-n^{-1}x)-1)$. The sequence $f_n(X)$ converges almost surely to $-X$. 
Furthermore, $|f_n(X)|$ is almost surely bounded from above by $|X|$ (because $|f_n(x)|\leq |x|$ if $\Re(x)\geq 0$) and 
$X$ is integrable (see Lemma \ref{lem2.1} below). Hence, using dominated convergence, we prove that $\bbE(f_n(X))$ converges to $-\bbE(X)$. This last fact is equivalent to equation (\ref{eq2.4}).\\
Finaly, equations (\ref{eq2.1}) and (\ref{eq2.4}) together yield 
\begin{equation}\label{eq2.5}
\lim_{n\rightarrow\infty}\bbE\left[\exp\left( i\sum_{j=1}^k \theta_j\Gamma_n(t_j)\right)\right]=\exp\left(-\bbE(X)\right).
\end{equation}
$\Box$\\
{\bf Remark :} The limit process is denoted by $\Gamma$. Its finite dimensional distributions are given by equation (\ref{eq1.12}), which is a consequence of equations (\ref{eq2.3}) and (\ref{eq2.5}) together. 

In the previous demonstration, we use the fact that $X$ is integrable. Since
$$|X|\leq (1+|\nu\tan(\frac{\pi\beta}{2})|)\sigma^\beta\left(\sum_{j=1}^k |\theta_j|\right)^\beta \int_{\bbsR} L_t(x)^{\beta}{\rm d}x$$
with $t=\max\{t_j\ ,\ 1\leq j\leq k\}$, the integrability of $X$ is a consequence of the following Lemma:

\begin{lem}\label{lem2.1}\ \\
For any $t\geq 0$, the random variable $\int_{\bbsR} L_t(x)^{\beta}{\rm d}x$ is integrable. 
\end{lem}
{\it Proof:\ }
In the case $\beta=2$, the following estimate from Khoshnevisan and Lewis (see \cite{LK98} cor 5.6 p.105)
$$\lim_{\lambda \rightarrow \infty}\lambda^{\alpha} \bbP\left(\int_{\bbsR} L_1^{2}(x){\rm d}x >\lambda\right)=c\in(0,\infty)$$
proves that the selfintersection local time of the $\alpha$-stable L\'evy process has bounded expectation at time $1$.
The self-similarity of the local time implies that 
$$\bbE\left( \int_{\bbsR} L_t^{2}(x){\rm d}x\right)=t^{2-\frac{1}{\alpha}}\bbE\left( \int_{\bbsR} L_1^{2}(x){\rm d}x\right)<\infty.$$
In the case $\beta\in(0,2)$, we use Hölder inequality,
\begin{eqnarray*}
\bbE\left(\int_{\bbsR} L_t^{\beta}(x){\rm d}x \right)
&=&\bbE\int_{\bbsR} L_t^{\beta}(x)1_{\{\sup_{0\leq s\leq t} |Y(s)|\geq |x| \}} {\rm d}x \\
&\leq& \left(\bbE\int_{\bbsR}  L_t^{2}(x){\rm d}x\right)^{\frac{\beta}{2}} \left( \int_{\bbsR} \bbP\left(\sup_{0\leq s\leq t} |Y(s)|\geq |x| \right){\rm d}x\right)^{1-\frac{\beta}{2}} 
\end{eqnarray*}
The last integral is bounded because the probability tail of suprema of $\alpha$-stable L\'evy process verifies
$$\lim_{\lambda \rightarrow \infty}\lambda^{\alpha}\bbP\left(\sup_{0\leq s\leq t} |Y(s)|\geq \lambda \right)=c'\in(0,\infty).$$
(see Theorem 10.51. p470 in \cite{ST94})
\\ $\Box$

\begin{pr}\label{pr1.2}
The sequence of process $\Gamma_n$ is tight in $\cC([0,\infty))$.
\end{pr} 
{\it Proof : } The case $\beta=2$ is straightforward: the process $\Delta$ is square integrable and for all $0\leq t_1<t_2$
\begin{eqnarray*}
\bbE\left[|\Delta(t_2)-\Delta(t_1)|^2\right]&=& \bbE\left[(\int_{\bbsR} (L_{t_2}(x)-L_{t_1}(x))\dd W_x)^2\right]\\ 
  &=& \sigma^2 \bbE\left[\int_{\bbsR} (L_{t_2}(x)-L_{t_1}(x))^2\dd x\right]\\ 
  &=& \sigma^2 \bbE\left[\int_{\bbsR} L_{t_2-t_1}^2(x)\dd x\right]\\ 
  &=& \sigma^2(t_2-t_1)^{2-\frac{1}{\alpha}}\bbE\left[\int_{\bbsR} L_{1}(x)^2\dd x\right]. 
\end{eqnarray*}
In the computation above, we first use  Itô formula to compute the second order moment of a stochastic integral with respect to Brownian motion. And then, we use of the self-similarity and the stationarity of the increments of the local time. 
Since $\Gamma_n$ is a sum of independent copies of $\Delta$,
\begin{eqnarray*}
\bbE\left[|\Gamma_n(t_2)-\Gamma_n(t_1)|^2\right]&=& \bbE\left[|n^{-\frac{1}{2}} \sum_{i=1}^n \Delta^{(i)}(t_2)-\Delta^{(i)}(t_1) |^2\right]\\ 
  &=&  \sigma^2(t_2-t_1)^{2-\frac{1}{\alpha}}\bbE\left[\int_{\bbsR} L_{1}(x)^2\dd x\right]. 
\end{eqnarray*}
Using Kolmogorov criterion, we deduce that the sequence $\Gamma_n$ is tight. 

In the case $0<\beta<2$, the proof has to be modified because the process $\Gamma_n$ has infinite variance. In order to prove the tightness of $\Gamma_n$, let us introduce a truncation method.
As a Lévy process, $W$ is a semimartingale whose Lévy-Itô decomposition is given by
\begin{equation}\label{eq2.6}
W_x=bx+\int_0^x\int_{|u|\leq 1} u(\mu-\bar\mu)(\dd u,\dd s)+\int_0^x\int_{|u|> 1} u\mu (\dd u,\dd s)
\end{equation}
where $\mu$ is a Poisson random measure on $\bbR\times\bbR$ with intensity $\bar\mu(\dd u,\dd x)=\lambda(\dd u)\otimes \dd x$,
$b$ is the drift and $\lambda$ stands for the stable Lévy measure on $\bbR$:
$$\lambda(\dd u)=\left(c_-1_{\{u<0\}}+c_+1_{\{u>0\}}\right) \frac{\dd u}{|u|^{\beta +1}},\ \ c_-,c_+\geq 0, c_-+c_+>0.$$
For some truncation level $R>1$, let $W^{(R^-)}$ and $W^{(R^+)}$ be the independent Lévy processes defined by 
$$W_x^{(R^-)}=\int_0^x\int_{|u|\leq R} u(\mu-\bar\mu)(\dd u,\dd s),\ \ \ W_x^{(R^+)}=\int_0^x\int_{|u|> R} u\mu (\dd u,\dd s).$$
The first one has a compactly supported Lévy measure and is a square-integrable martingale with infinitely many jumps bounded by $R$ on each compact time interval, whereas the second one is a compound Poisson process. The Lévy-Itô decompostion (\ref{eq2.6}) rewrites as 
$$W_x=b_R x+W_x^{(R^-)}+W_x^{(R^+)}$$
where $b_R=b+\int_{1<|y|\leq R} u\lambda(\dd u)$ is a drift depending on $R$.

This decomposition of the stable scenery yields the following decomposition of the stable process in random scenery:
$$\Delta(t)=b_Rt + \Delta^{(R^-)}(t) + \Delta^{(R^+)}(t),$$
with 
$$\Delta^{(R^-)}(t)=\int_{\bbsR} L_t(x)\dd W_x^{(R^-)}, \ \ \ \Delta^{(R^+)}(t)=\int_{\bbsR} L_t(x)\dd W_x^{(R^+)} .$$
The process $\Delta^{(R^-)}$ is square integrable and for any $0\leq t_1<t_2$,
\begin{eqnarray}
\bbE\left[(\Delta^{(R^-)}(t_2)- \Delta^{(R^-)}(t_1))^2\right]&=&\bbE\left[(\int_{\bbsR} (L_{t_2}(x)-L_{t_1}(x))\dd W_x^{(R^-)})^2\right] \nonumber \\
 &=& \bbE\left[\int_{\bbsR} (L_{t_2}(x)-L_{t_1}(x))^2\dd <W^{(R^-)},W^{(R^-)}>_x\right] \nonumber \\
 &=&  \frac{c_-+c_+}{2-\beta} R^{2-\beta} \bbE\left[\int_{\bbsR} L_{t_2-t_1}(x)^2\dd x\right] \nonumber \\
 &=&  \frac{c_-+c_+}{2-\beta} R^{2-\beta} (t_2-t_1)^{2-\frac{1}{\alpha}}\bbE\left[\int_{\bbsR} L_{1}(x)^2\dd x\right] \label{eq2.7}
\end{eqnarray}
Here we have used the following expression of the bracket angle
$$  <W^{(R^-)},W^{(R^-)}>_x=\int_0^x \int_{|u|\leq R} u^2 \bar\mu(\dd u, \dd s)=\frac{c_-+c_+}{2-\beta} R^{2-\beta} x,$$
 and also the properties of stationarity and self-similarity of the stable process that goes through its local time.\\
 
On the other hand, since $W_x^{(R^+)}$ is a pure jump process, the stochastic integral defining $\Delta^{(R^+)}$ is constructed in the Lebesgue-Stieljes sense and is equal to
$$\Delta^{(R^+)}(t)=\int_{\bbsR} \int_{|u|> R} uL_t(s)\mu (\dd u,\dd s).$$
Let $T\geq 0$. Suppose that the Lévy-process $W$ has no jump of size greater than $R$ on the (bounded) support of $L_T$, then the process 
$\Delta^{(R^+)}$ is identically $0$ on $[0,T]$. Hence,
$$\bbP\left( \Delta^{(R^+)} \equiv 0 {\rm\ on \ } [0,T] \right) \geq \bbP\left(\mu \{(u,x)\in\bbR^2 |\ |u|\geq R {\rm\ and\ } L_T(x)\neq 0\}=0\right) .$$
Conditionally on the support of the local time, the random variable  
$$\mu \{(u,x)\in\bbR^2 |\ |u|\geq R {\rm\ and\ } L_T(x)\neq 0\}$$
has Poisson distribution with parameter 
\begin{eqnarray*}
\bar\mu\left(\supp(L_T)\times(\bbR\setminus [-R,R])\right)&=&
\int_{\bbsR}1_{\{L_T(x)\neq 0\}} \dd x \ \int_{|u|>R} \lambda(y)\dd u \\
&\leq&  2\sup_{0\leq t\leq T} |Y(t)| \frac{c^++c^-}{\beta}R^{-\beta}.
\end{eqnarray*}
Hence, the probability that $\Delta^{(R^+)}$ is identically $0$ on $[0,T]$ is underestimated by
\begin{eqnarray}
\bbP\left( \Delta^{(R^+)} \equiv 0 {\rm\ on \ } [0,T] \right) &\geq& \bbE\left[\exp\left(-2\frac{c^++c^-}{\beta}R^{-\beta}\sup_{0\leq t\leq T} |Y(t)| \right)\right] \nonumber \\
&\geq& 1-2\frac{c^++c^-}{\beta}R^{-\beta} \bbE\left(\sup_{0\leq t\leq T} |Y(t)| \right) \label{eq2.8}
\end{eqnarray}

We now consider a corresponding decomposition for $\Gamma_n$. Let $R>1$ and introduce the truncation level $R_n=Rn^{\frac{1}{\beta}}$.
Each random walk in random scenery can be decomposed as above, what yields,  
$$\Delta^{(i)}(t)=b_{R_n}t + \Delta^{(i,R_n^-)}(t) + \Delta^{(i,R_n^+)}(t).$$
The process $\Gamma_n$ writes
\begin{equation}\label{eq2.9}
\Gamma_n(t)=n^{1-\frac{1}{\beta}} b_{R_n} t +\Gamma_n^{(R^-)}(t)+\Gamma_n{(R^+)}(t),
\end{equation}
with 
$$\Gamma_n^{(R^-)}(t)=n^{-\frac{1}{\beta}}\sum_{i=1}^n \Delta^{(i,R_n^-)}(t), \ \ \Gamma_n^{(R^+)}(t)=n^{-\frac{1}{\beta}}\sum_{i=1}^n \Delta^{(i,R_n^+)}(t)$$
The tightness of the sequence $\Gamma_n$ is a consequence of equation (\ref{eq2.9}) and of the following assertions :\\
$A_1$ - the sequence  $n^{1-\frac{1}{\beta}} b_{R_n} ,\ n\geq 1$ is bounded,\\
$A_2$ - for each $R>1$, the sequence of process $\Gamma_n^{(R^-)}, n\geq 1$ is tight, \\
$A_3$ - for each $T>0$, the probability that $\Gamma_n^{(R^+)}\equiv 0 $ on $[0,T]$ satisfies
$$\lim_{R\rightarrow\infty}\limsup_{n\rightarrow\infty} \bbP\left(\Gamma_n^{(R^+)}\equiv 0 {\rm\ on\ }[0,T]\right) =1.$$ 
Proof of assertion $A_1$: we compute
$$b_R=b+\int_{1<|y|\leq R} u\lambda(\dd u)=b+\frac{c^++c^-}{1-\beta}(R^{1-\beta}-1),$$
and hence
$$n^{1-\frac{1}{\beta}} b_{R_n}=(b-\frac{c^++c^-}{1-\beta})n^{1-\frac{1}{\beta}}+\frac{c^++c^-}{1-\beta}R^{1-\beta}.$$
In the case $\beta\leq 1$, this quantity is trivially bounded. In the case $\beta>1$, this is still true because
$$\bbE(W_1)=b-\frac{c^++c^-}{1-\beta}=0.$$

Proof of assertion $A_2$: Let $R>1$. For any $0<t_1<t_2$, equation (\ref{eq2.7}) implies 
\begin{eqnarray}
\bbE\left[(\Gamma_n^{(R^-)}(t_2)- \Gamma_n^{(R^-)}(t_1))^2\right]&=& n^{1-\frac{2}{\beta}}\bbE\left[(\Delta^{(R_n^-)}(t_2)- \Delta^{(R_n^-)}(t_1))^2\right] \nonumber\\
 &=&  n^{1-\frac{2}{\beta}} \frac{c_-+c_+}{2-\beta} R_n^{2-\beta} (t_2-t_1)^{2-\frac{1}{\alpha}}\bbE\left[\int_{\bbsR} L_{1}(x)^2\dd x\right]\nonumber\\
 &=&  \frac{c_-+c_+}{2-\beta} R^{2-\beta} (t_2-t_1)^{2-\frac{1}{\alpha}}\bbE\left[\int_{\bbsR} L_{1}(x)^2\dd x\right]\label{eq2.10}.
\end{eqnarray}
Using Kolmogorov criterion, this estimate implies that $\Gamma_n^{(R^-)}$ is tight.\\
Proof of assertion $A_3$: 
As $\Gamma_n^{(R^+)}$ is a sum of $n$ independent copies of $\Delta^{(R_n^+)}$, it vanishes on $[0,T]$ as soon as each summand does. Hence, using equation (\ref{eq2.8}), we get
\begin{eqnarray*}
\bbP\left( \Gamma_n^{(R^+)}\equiv 0 {\rm\ on\ }[0,T] \right)&\geq& \left[\bbP\left( \Delta^{(R_n^+)}\equiv 0 {\rm\ on\ }[0,T] \right)\right]^n \\
 &\geq& \left[  1-2\frac{c^++c^-}{\beta}R_n^{-\beta} \bbE\left(\sup_{0\leq t\leq T} |Y(t)| \right)\right]^n.
\end{eqnarray*}
Take into account the scaling $R_n=Rn^{\frac{1}{\beta}}$, and let $n$ go to infinity:
$$\limsup_{n\rightarrow\infty} \bbP\left(\Gamma_n^{(R^+)}\equiv 0 {\rm\ on\ }[0,T]\right)\geq \exp\left[-2\frac{c^++c^-}{\beta}R^{-\beta} \bbE\left(\sup_{0\leq t\leq T} |Y(t)| \right)\right].$$
This last estimate implies the asymptotic given in $A_3$ as $R$ goes to infinity.
$\Box$\\

{\bf Remark: } We have proven the tightness of the sequence $\Gamma_n,n\geq 1$ by using Kolmogorov's criterion. As a by product of  estimation (\ref{eq2.10}), we can deduce that the limit $\Gamma$ is locally $\gamma$-H\"older for all $\gamma <\frac{\alpha-1}{2\alpha}  $. However, this is not sharp
and Theorem \ref{theo5} gives a better result about the regularity of the trajectories of the process $\Gamma$.

\section{The random rewards schema \\ Proof of Theorem \ref{theo2}}
\subsection{Preliminaries : local times, number of self-intersections, and range.}
Let $x\in\bbZ$ and $n\geq 1$. The local time $N_n(x)$ of the random walk $(S_k)_{k\geq 0}$ at point $x$ up to time $n$ is defined by
$$N_n(x)=\sum_{k=0}^n 1_{\{S_k=x\}}.$$
It represents the amount of time the walk spends at state $x$ up to time $n$. We extend this definition to noninteger
time $s\geq 0$ by linear interpolation:
$$N_s(x)=N_{[s]}(x)+(s-[s])(N_{[s]+1}(x)-N_{[s]}(x).$$
The random walk in random scenery defined by eqs (\ref{eq1.06}) and (\ref{eq1.07}) writes for all $s\geq 0$  
\begin{equation}\label{eq3.1}
Z_s=\sum_{x\in Z} N_s(x)\xi_x 
\end{equation}
where the collection of random variables $\{N_s(x), x\in\bbZ\}$ and $\{\xi_x, x\in\bbZ\}$ are independent.
The number of self-intersections $V_n$ of the random walk up to time $n$ is defined by
\begin{equation}\label{eq3.2}
V_n=\sum_{0\leq i,j\leq n} 1_{\{S_i=S_j\}}=\sum_{x\in\bbsZ}N_n(x)^2.
\end{equation}
The range $R_n$ of the random walk up to time $n$ is defined by
\begin{equation}\label{eq3.3}
R_n=\card\left\{S_k;0\leq k\leq n\right\}=\sum_{x\in\bbsZ}1_{\{N_n(x)\neq 0\}}.
\end{equation}
It represents the number of different states visited by the random walk up to time $n$.\\
The definitions (\ref{eq3.2}) and (\ref{eq3.3}) extend obviously to noninteger time $s\geq 0$.

The following Lemma yields some estimations of these quantities.
\begin{lem}\label{lem3.1}\ \\
\begin{itemize}
\item The following convergence in probability holds 
\begin{equation}\label{eq3.4}
n^{-\delta}\sup_{x\in\bbsZ} N_n(x) \mathop{\longrightarrow}_{n\rightarrow\infty}^{{\rm P}} 0.
\end{equation} 
\item  For any $p\in [1,+\infty)$, there exists some constant $C$ such that for all $n\geq 1$, 
\begin{equation}\label{eq3.5}
\bbE\left(V_n^p\right)\leq C n^{(2-\frac{1}{\alpha})p}.
\end{equation}
\item For any $p\in [1,\alpha)$, there exists some constant $C$ such that for all $n\geq 1$
\begin{equation}\label{eq3.6}
\bbE\left( R_n^p\right)\leq Cn^{\frac{p}{\alpha}}.
\end{equation}
\end{itemize}
\end{lem}
\noindent
{\it Proof:}\ \\
$\bullet$ The proof is given in \cite{KS79} Lemma 4 p.12.\\
$\bullet$ Since if $0<p'<p$, $\bbE(|X|^{p'})\leq \bbE(|X|^p)+1$, we can suppose that $p$ is an integer $\geq 1$. The number of self-intersection up to time $n$ is bounded from above by
$$V_n\leq 2\sum_{0\leq i\leq j\leq n} 1_{\{S_i=S_j\}}.$$
Using Minkowsky inequality,
\begin{equation}\label{eq3.7}
|\!|V_n|\!|_p \leq 2 \sum_{i=0}^n |\!| \sum_{j=i}^n I(S_i=S_j)|\!|_p,
\end{equation}
where $|\!|X|\!|_p=\bbE(|X|^p)^{1/p}$.
For fixed $i$, the stationarity of the random walk's increments implies that the distribution of $\sum_{j=i}^n 1_{\{S_i=S_j\}}$ and 
$\sum_{j=0}^{n-i} 1_{\{S_i=0\}}=N_{n-i}(0)$ are equal. Since $N_{n-i}(0)\leq N_n(0)$, equation (\ref{eq3.7}) yields
\begin{equation}\label{eq3.8}
\bbE(V_n^p)\leq 2^p n^p \bbE\left(N_n(0)^p\right).
\end{equation}
Lemma 1 p.12 of \cite{KS79} states that there exists some $C>0$ such that for all $n\geq 1$
\begin{equation}\label{eq3.9}
 \bbE\left(N_{n}(0)^p\right)\leq C n^{p(1-\frac{1}{\alpha})}
\end{equation}
Equations (\ref{eq3.8}) and (\ref{eq3.9}) together yield equation (\ref{eq3.5}).\\
$\bullet$ We suppose $p>1$, the case $p=1$ follows. For $n\geq 0$, define the random variable $M_n$ by 
\begin{equation}\label{eq3.10}
M_n=\max_{0\leq k\leq n} |S_k|.
\end{equation}
Up to time $n$, the random walk $(S_k)_{k\geq 0}$ remains in the set $\left\{x\in\bbZ\ ;\ |x|\leq M_n\right\}$, and hence
\begin{equation}\label{eq3.11}
R_n\leq 2M_n+1\ \ {\rm a.s.}.
\end{equation}
Since $\left(|S_n|\right)_{n\geq 0}$ is a non-negative submartingale and since $\bbE(|S_n|^p)<\infty$ if $p\in (1,\alpha)$, we can apply the maximal inequality (Theorem 20 p 11 of \cite{Pr04}). This yields that for any $p\in(1,\alpha)$,
\begin{equation}\label{eq3.12}
\bbE\left(M_n^p\right)\leq \left(\frac{p}{p-1}\right)^p\bbE\left(|S_n|^p\right).
\end{equation}
We need the following estimate (equation 5s p.684 in \cite{LGR91}): for any $p\in(1,\alpha)$, there exists some constant $C>0$ such that for all $n\geq 1$,
\begin{equation}\label{eq3.13}
\bbE\left(|S_n|^p\right)\leq Cn^{\frac{p}{\alpha}}.
\end{equation}
Combining equations (\ref{eq3.11}),(\ref{eq3.12}) and (\ref{eq3.13}), we prove equation (\ref{eq3.6}) for $p>1$. 
$\Box$

We recall the following lemma from Kesten and Spitzer which is in the bulk of the proof. 
\begin{lem}\label{KS}
For all $(\theta_1,\cdots,\theta_k)\in\bbR^k$, $(t_1,\cdots,t_k)\in [0,+\infty)^k$, $\sigma>0$, $\beta\in(0,2]$ and $\nu\in[-1,1]$ (with $\nu=0$ when $\beta=1$), the distribution of
$$X_n=\sigma^{\beta} n^{-\delta\beta}\sum_{x\in\bbsZ}\big|\sum_{j=1}^k \theta_j N_{nt_j}(x)\big|^\beta \left(1- i\nu \tan(\frac{\pi\beta}{2})\sgn(\sum_{j=1}^k \theta_j N_{nt_j}(x))\right) $$
converges weakly as $n\rightarrow\infty$ to $X$ defined by equation (\ref{eq2.3}).
\end{lem}

\begin{lem}\label{KSbis}
Furthermore, for all $p\in [1,\alpha(1-\frac{\beta}{2})^{-1})$  (if $\beta=2$, read for all $p\in[1,\infty)$),
there exists some constant $C>0$ such that for all $n\geq 1$,
$$\bbE(|X_n|^p)\leq C .$$
\end{lem}
\noindent
{\it Proof:}\ \\ 
Let $T=\max(t_1,\cdots,t_n)$ and $\Theta=\sum_{j=1}^k |\theta_j|$. 
The random variables $|X_n|$ is bounded above by 
$$(1+|\nu\tan(\frac{\pi\beta}{2})|)\Theta^\beta n^{-\delta\beta} \sum_{x\in\bbsZ} N_{[nT]+1}^\beta(x).$$
In the case $\beta=2$, this quantitiy is equal to $\Theta^2 n^{\frac{1}{\alpha}-2} V_{[nT]+1}$, and in this case, Lemma \ref{KSbis} is a consequence of equation (\ref{eq3.5}).\\
In the case $\beta<2$,  Hölder inequality yields
$$\sum_{x\in\bbsZ} N_{[nT]+1}^\beta(x)\leq  \left(\sum_{x\in\bbsZ} 1_{\{N_{[nT]+1}(x)\neq0\}}\right) ^{1-\frac{\beta}{2}}\left(\sum_{x\in\bbsZ} N_{[nT]+1}^2(x)\right)^{\frac{\beta}{2}} =R_{[nT]+1}^{1-\frac{\beta}{2}}V_{[nT]+1}^{\frac{\beta}{2}}.$$
Hence, up to a multiplicative constant, the expectation $\bbE(|X_n|^p)$ is overestimated by
$$\bbE\left[ \left(n^{-\delta\beta} \sum_{x\in\bbsZ} N_{[nT]+1}^\beta(x)\right)^p\right]\leq 
\bbE \left[ \left(n^{-\frac{1}{\alpha}}R_{[nT]+1} \right)^{p(1-\frac{\beta}{2})}\left(n^{\frac{1}{\alpha}-2}V_{[nT]+1}\right)^{p\frac{\beta}{2}}\right]$$
We now apply Hölder inequality another time, where the random variables are positive and the expectations eventually $+\infty$: for all $\gamma>\frac{p\beta}{2}$, 
\begin{eqnarray*}
& &\bbE \left[ \left(n^{-\frac{1}{\alpha}}R_{[nT]+1} \right)^{p(1-\frac{\beta}{2})}\left(n^{\frac{1}{\alpha}-2}V_{[nT]+1}\right)^{p\frac{\beta}{2}}\right]\\
&\leq& \bbE \left[ \left(n^{-\frac{1}{\alpha}}R_{[nT]+1} \right)^{p(1-\frac{\beta}{2})\frac{2\gamma}{2\gamma-p\beta}}\right]^{\frac{2\gamma-p\beta}{2\gamma}} \bbE\left[\left(n^{\frac{1}{\alpha}-2}V_{[nT]+1}\right)^{\gamma}\right]^{\frac{p\beta}{2\gamma}}
\end{eqnarray*}
We can choose $\gamma$ such that this last term is bounded as $n\rightarrow\infty$.
>From Lemma \ref{lem3.1}, this is the case if the condition $p(1-\frac{\beta}{2})\frac{2\gamma}{2\gamma-p\beta}<\alpha$ is satisfied. Since $ p(1-\frac{\beta}{2})<\alpha$, this last condition is satisfied for $\gamma$ large enough.
$\Box$

\subsection{Convergence of the finite-dimensional distributions.}
We define the process $G_n$ by 
\begin{equation}\label{eq3.14}
G_n(t)= c_n^{-\frac{1}{\beta}} \sum_{i=1}^{c_n} D_n^{(i)}(t),\ \ t\geq 0,
\end{equation}
where $D_n^{(i)}$ is the i-th random walk in random scenery properly rescaled and defined by (\ref{eq1.10}). Theorem \ref{theo2} states that $G_n$ converges weakly to $\Gamma$ in $\cC([0,\infty))$.\\

The convergence of the finite-dimensional marginals is proven by studying the asymptotic behaviour of their characteristic function. \\
Let $\lambda$ be the characteristic function of the variables $\xi_k^{(i)}$ defined by
$$\lambda(u)=\bbE\left(\exp(iu\xi_1^{(1)})\right).$$
Since the random variables $\xi_k^{(i)}$ are in the domain of attraction of $Z_{\beta}$, 
\begin{equation}\label{eq3.15}
\lambda(u)=\bar\lambda (u)+o(|u|^\beta)\ \ ,\ \ {\rm as}\ u\rightarrow 0,
\end{equation}
where $\bar\lambda$ is the characteristic function of $Z_beta$ given by equation (\ref{eq1.02}).

\begin{pr}\label{pr1}
The finite dimensional distributions of $(G_n(t))_{t\geq 0}$ converge weakly as $n\rightarrow\infty$ to those of $(\Gamma(t))_{t\geq 0}$ defined in equation (\ref{eq1.09}).
\end{pr}

{\bf Remark } It is worth noting that the properties of stability, selfsimilarity and stationarity can be seen in the characteristic function of the marginals.

For sake of clarity, we divide the proof of Proposition \ref{pr1} into three lemmas.
In Lemma \ref{lem1} below, we express the characteristic function in terms of the local time of the random walks $S_n$ and of the characteristic function $\lambda$ of the scenery. In Lemma \ref{lem2}, we estimate the error committed when replacing the scenery by a stable scenery. Finally, in Lemma \ref{lem3}, we prove the proposition in the case of a stable scenery.

\begin{lem}\label{lem1}
For all $(\theta_1,\cdots,\theta_k)\in \bbR^k$ , $(t_1,\cdots,t_k)\in [0,+\infty)^k$,
$$ \bbE\left[\exp\left(i\sum_{j=1}^k \theta_j G_n(t_j)\right) \right]
= \left(\bbE \left[ \prod_{x\in\bbsZ} \lambda \left( c_n^{-\frac{1}{\beta}}n^{-\delta} \sum_{j=1}^k \theta_j N_{nt_j}(x)\right)\right]\right )^{c_n}$$
\end{lem}
\noindent
{\it Proof:}\ \ 
From equation (\ref{eq3.14}) and the independence of the processes $\left(D_n^{(i)}(t)\right)_{t\geq 0}$, 
\begin{equation}\label{eq3.15b}
\bbE\left[\exp\left(i\sum_{j=1}^k \theta_j G_n(t_j)\right) \right]
=\bbE\left[\exp\left(ic_n^{-\frac{1}{\beta}}\sum_{j=1}^k \theta_j D_{n}^{(1)} (t_j)\right) \right]^{c_n}
\end{equation}
Hence we have to compute the characteristic function of a single random walk in random 
scenery. In the sequel, we omit the superscript $(1)$, and write $D_n(t)$ (resp. $S_k$, $\xi_k\ \cdots$)
instead of $D_n^{(1)}(t)$ (resp. $S_k^{(1)}$, $\xi_k^{(1)}\ \cdots$). 
Using equation (\ref{eq1.10}) and (\ref{eq3.1}), the characteristic function writes
$$\bbE\left[\exp\left(ic_n^{-\frac{1}{\beta}}\sum_{j=1}^k \theta_j D_{n}(t_j)\right) \right]=
\bbE\left[\prod_{x\in \bbsZ} \exp\left(ic_n^{-\frac{1}{\beta}}n^{-\delta}\sum_{j=1}^k \theta_j N_{nt_j}(x)\xi_x\right) \right].$$
Using the independence of the random walk $(S_n)_{n\geq 0}$ and the random scenery $(\xi_x)_{x\in\bbsZ}$, and performing integration on the scenery yields
\begin{equation}\label{eq3.16}
\bbE\left[\exp\left(ic_n^{-\frac{1}{\beta}}\sum_{j=1}^k \theta_j D_{n}(t_j)\right) \right]=\bbE\left[\prod_{x\in \bbsZ} \lambda\left(c_n^{-\frac{1}{\beta}}n^{-\delta}\sum_{j=1}^k \theta_j N_{nt_j}(x)\right) \right].
\end{equation}
Equations (\ref{eq3.15b}) and (\ref{eq3.16}) together yield Lemma \ref{lem1}
$\Box$

\begin{lem}\label{lem2}
For all $(\theta_1,\cdots,\theta_k)\in \bbR^k$ , $(t_1,\cdots,t_k)\in [0,+\infty)^k$,
$$\bbE\left[ \prod_{x\in\bbsZ} \lambda \left( c_n^{-\frac{1}{\beta}}n^{-\delta} \sum_{j=1}^k \theta_j N_{nt_j}(x)\right)\right]
=\bbE\left[ \prod_{x\in\bbsZ} \bar\lambda \left( c_n^{-\frac{1}{\beta}}n^{-\delta} \sum_{j=1}^k \theta_j N_{nt_j}(x)\right)\right]+o(c_n^{-1})$$
\end{lem}
\noindent
{\it Proof:}\ \\
Denote by $U_{n}(x)$  the random variables defined by
$$U_{n}(x)=n^{-\delta} \sum_{j=1}^k \theta_j N_{nt_j}(x)\ \ ,\ \ x\in\bbZ$$
With these notations, Lemma \ref{lem2} is equivalent to
$$\lim_{n\rightarrow\infty} \bbE\left[ c_n\left( \prod_{x\in\bbsZ} \lambda \left( c_n^{-\frac{1}{\beta}}U_{n}(x)\right)- \prod_{x\in\bbsZ} \bar \lambda \left( c_n^{-\frac{1}{\beta}}U_{n}(x)\right) \right) \right]= 0. $$
Note that the products although indexed by $x\in\bbZ$ have only a finite number of factors different from $1$.
And furthermore, all factors are complex numbers in $\bar{\bbD} = \left\{z\in\bbC | \ |z|\leq 1\right\}$.
We use the following inequality : let $(z_i)_{i\in I}$ and $(z_i')_{i\in I}$ two families of complex numbers in $\bar{\bbD}$ such that all terms are equal to one, except a finite number of them. Then
$$\left|\prod_{i\in I}z'_i -\prod_{i\in I}z_i\right| \leq \sum_{i\in I} \left|z'_i - z_i\right|.$$
This yields 
\begin{eqnarray}\label{eq3.17}
 & &\left|\ \prod_{x\in\bbsZ} \lambda \left( c_n^{-\frac{1}{\beta}}U_{n}(x)\right)- \prod_{x\in\bbsZ} \bar \lambda \left( c_n^{-\frac{1}{\beta}}U_{n}(x)\right)\ \right| \nonumber \\
 &\leq& \  \sum_{x\in\bbsZ}  \ \left|\ \lambda \left( c_n^{-\frac{1}{\beta}}U_{n}(x)\right)- \bar\lambda \left( c_n^{-\frac{1}{\beta}}U_{n}(x)\right)\  \right|.
\end{eqnarray}
We define a continuous and bounded function $g$ by $g(0)=0$ and
$$g(v)=|v|^{-\beta}\left|\lambda(v)-\bar\lambda(v)\right|\ \ ,\ \ v\neq 0,$$
so that for any $x\in\bbZ$,
$$\left| \lambda \left( c_n^{-\frac{1}{\beta}}U_{n}(x)\right)-  \bar \lambda \left( c_n^{-\frac{1}{\beta}}U_{n}(x)\right) \right| = c_n^{-1}|U_{n}(x)|^\beta g(c_n^{-\frac{1}{\beta}}U_{n}(x)).$$
In order to obtain uniform estimations, define the random variable $U_{n}$ by
$$U_{n}=\max_{x\in\bbsZ} |U_{n}(x)|,$$
and the function $\tilde g:[0,+\infty)\rightarrow [0,+\infty)$ by
$$\tilde g(u) =\sup_{|v|\leq u} |g(v)|.$$
Note that $\tilde g$ is continuous, bounded and vanishes at $0$.
Then, for any $x\in \bbZ$,
\begin{equation}\label{eq3.18}
\left| \lambda \left( c_n^{-\frac{1}{\beta}}U_{n}(x)\right)-  \bar \lambda \left( c_n^{-\frac{1}{\beta}}U_{n}(x)\right) \right|\leq c_n^{-1}|U_{n}(x)|^\beta \tilde g(c_n^{-\frac{1}{\beta}}U_{n}).
\end{equation}
Equations (\ref{eq3.17}) and (\ref{eq3.18}) together yield
$$c_n\left| \prod_{x\in\bbsZ} \lambda \left( c_n^{-\frac{1}{\beta}}U_{n}(x)\right)- \prod_{x\in\bbsZ} \bar \lambda \left( c_n^{-\frac{1}{\beta}}U_{n}(x)\right) \right|\leq \tilde g(c_n^{-\frac{1}{\beta}}U_{n})\sum_{x\in\bbsZ}|U_{n}(x)|^\beta.$$
>From Lemma \ref{lem3.1}, $U_{n}$ converge in probability to $0$ as $n\rightarrow\infty$. Since $\tilde g$ is continuous and vanishes at $0$, $\tilde g(c_n^{-\frac{1}{\beta}}U_{n})$ converges also in probability to $0$. From Lemma \ref{KS} (with $\nu=0$),
$\sum_{x\in\bbsZ}|U_{n}(x)|^\beta$ converges in distribution. Hence, the product  
$$\tilde g(c_n^{-\frac{1}{\beta}}U_{n}) \sum_{x\in\bbsZ}|U_{n}(x)|^\beta$$
converges to $0$ in probability.
We finally prove that the expectation converges also to $0$. We use uniform integrability and prove that the sequence of random variable $\tilde g(c_n^{-\frac{1}{\beta}}U_{n}) \sum_{x\in\bbsZ}|U_{n}(x)|^\beta$ is bounded in $L^p$ for some $p>1$. The function $\tilde g$ is bounded on $[0,+\infty)$, and Lemma \ref{KSbis} states that $\sum_{x\in\bbsZ}|U_{n}(x)|^\beta$ is bounded in $L^p$ for some $p>1$.
$\Box$

\begin{lem}\label{lem3}
For all $(\theta_1,\cdots,\theta_k)\in \bbR^k$ , $(t_1,\cdots,t_k)\in [0,+\infty)^k$,
$$\bbE\left[ \prod_{x\in\bbsZ} \bar\lambda \left( c_n^{-\frac{1}{\beta}}n^{-\delta} \sum_{j=1}^k \theta_j N_{
nt_j}(x)\right)\right] =1- c_n^{-1}\ \bbE(X) +o(c_n^{-1}),$$
where $X$ is defined in Lemma \ref{KS}.
\end{lem}
\noindent
{\it Proof:}\ \\
Recall the definition of the random variable $X_n$ from Lemma \ref{KS} and of the characteristic function $\bar\lambda$ from equation (\ref{eq1.02}). With this notations, Lemma \ref{lem3} is equivalent to 
$$\lim_{n\rightarrow +\infty}  \bbE\left(f_n(X_n) \right)=\bbE(X),$$
where $f_n$ is the function defined on $\bbC$ by 
$$f_n(x)=c_n\left(1-\exp(-c_n^{-1}x)  \right).$$
It is easy to verify that the sequence of functions $f_n$ satisfies the following property: for every $x$, for every sequence $(x_n)_{n\geq 1}$ converging to $x$, 
$$\lim_{n\rightarrow\infty} f_n(x_n)=x.$$
Furthermore Lemma \ref{KS} states that the sequence $(X_n)_{n\geq 1}$ converges in distribution to $X$ when $n\rightarrow\infty$. Using the continuous mapping Theorem (Therorem 5.5 of \cite{Bi68}), we prove the weak convergence of the sequence of random variables $f_n(X_n)$ to $X$. \\
In order to prove the convergence of the expectation, we use uniform integrability and prove that there exists some $p>1$ such that  $\bbE\left(\left|f_n(X_n)\right|^p\right)$ is bounded. The integral representation
$$f_n(x)=x\int_0^1 \exp(-uc_n^{-1}x){\rm d}u$$
shows that for $\Re(x)\geq 0$, 
$$|f_n(x)|\leq |x|.$$
Hence Lemma \ref{KSbis} implies that $\left|f_n(X_n)\right|\leq |X_n|$ is bounded in $L^p$ for some $p>1$ and this proves the Lemma. $\Box$

\noindent
{\it Proof of Proposition \ref{pr1}:}\ \\
Using Lemmas \ref{lem1}, \ref{lem2} and \ref{lem3}, we have
\begin{eqnarray*}
\bbE\left[\exp\left(i\sum_{j=1}^k \theta_j G_n(t_j)\right) \right]
&=& \left(\bbE \left[ \prod_{x\in\bbsZ} \lambda \left( c_n^{-\frac{1}{\beta}}n^{-\delta} \sum_{j=1}^k \theta_j N_{nt_j}(x)\right)\right]\right )^{c_n}\\
&=&\left(\bbE \left[ \prod_{x\in\bbsZ} \bar\lambda \left( c_n^{-\frac{1}{\beta}}n^{-\delta} \sum_{j=1}^k \theta_j N_{nt_j}(x)\right)\right] +o(c_n^{-1})\right)^{c_n}\\
&=&\left(1- c_n^{-1}\ \bbE(X) +o(c_n^{-1})\right)^{c_n}.
\end{eqnarray*}
This last quantity has limit $\exp\left(-\bbE(X)\right)$ as $n\rightarrow \infty$.
This proves Propostition \ref{pr1}. $\Box$

\subsection{Tightness}
The weak convergence stated in Therorem \ref{theo1} is a consequence of Proposition \ref{pr1} and the following proposition :
\begin{pr}\label{pr2}
The family of processes $\left(\Gamma_{n}(t)\right)_{t\geq 0}$ is tight in $\cC([0,\infty))$.
\end{pr}

The tightness will be proven by approximating $G_n$ by processes with finite variance obtained with truncated random sceneries. Let $a$ be a positive real. We decompose the scenery $(\xi_x^{(i)})_{x\in\bbsZ,i\geq 1}$ into two parts
$$\xi_x^{(i)}=\bar\xi_{a,x}^{(i)}+\hat\xi_{a,x}^{(i)},$$
where $(\bar\xi_{a,x}^{(i)})$ denote the $i$-th truncated scenery defined by
$$\bar\xi_{a,x}^{(i)}=\xi_x^{(i)} 1_{\{ |\xi_x^{(i)}|\leq a\}},$$
and $\hat\xi_{a,x}^{(i)}$ the remainder scenery  
$$\hat\xi_{a,x}^{(i)}=\xi_x^{(i)} 1_{\{ |\xi_x^{(i)}| > a\}}.$$

In the sequel, we will use intensively the following estimates that are adapted from Feller \cite{Fe71}.
\begin{lem}\label{lem4}\ \\
\begin{itemize}
\item there exists $C$ such that for any  $a>0$,
$$|\bbE\left(\bar\xi_{a,x}^{(i)}\right)|\leq Ca^{1-\beta}$$
$$ \bbE\left(|\bar\xi_{a,x}^{(i)}|^2\right)\leq Ca^{2-\beta}.$$
\item there exists $C$ such that 
$$\bbP\left(\hat \xi_{a,x}^{(i)} \neq 0 \right) \leq  C a^{-\beta}.$$
\end{itemize}
\end{lem}
\noindent
{\it Proof:}\ \\
These estimations rely on the behaviour of the tail of the random variable $\xi_x^{(i)}$: there exist some constant $A_1$ and $A_2$ such that
$$\lim_{u\rightarrow +\infty} u^\beta\bbP(\xi_x^{(i)}>u)= A_1,$$
$$ \lim_{u\rightarrow +\infty} u^\beta\bbP(\xi_x^{(i)}<-u)= A_2.$$
We prove the lemma in the case $0<\beta<2$. In this case, we have $A_1\neq 0$ or $A_2\neq 0$. 
In the case $\beta=2$,  the constants $A_1$ and $A_2$ are zero and the demonstration is easily modified.\\
$\bullet$ The estimation of $ \bbE(\bar\xi_{a,x}^{(i)})$ comes from the formula
$$ \bbE(\bar\xi_{a,x}^{(i)})=\int_0^a \left(\bbP(\xi_x^{(i)}1_{\{|\xi_x^{(i)}|\leq a\}}>u)-\bbP(\xi_x^{(i)}1_{\{|\xi_x^{(i)}|\leq a\}}<-u)\right){\rm d}u.$$
In the case $\beta<1$, $| \bbE(\bar\xi_{a,x}^{(i)})|$ is bounded from above by  
$$\int_0^{a} \bbP(|\xi_x^{(i)}|>u){\rm d}u$$ 
and this integral is equivalent as $a\rightarrow +\infty$ to $Ca^{1-\beta}$.\\
In the case $\beta>1$,  $\xi_x^{(i)}$ has expectation $0$. This implies that $| \bbE(\bar\xi_{a,x}^{(i)})|$ is bounded from above by
$$\int_a^{\infty} \bbP(|\xi_x^{(i)}|>u){\rm d}u$$
and this integral is equivalent as $a\rightarrow +\infty$ to $Ca^{1-\beta}$.\\
In the case $\beta=1$, the symmetry condition (\ref{eq1.03}) means exactly that $| \bbE(\bar\xi_{a,x}^{(i)})|$ is bounded from above by some constant. \\
The second moment of the truncated scenery is overestimated by
\begin{eqnarray*}
\bbE\left(|\bar\xi_{a,x}^{(i)}|^2\right)&\leq& \bbE\left(|\xi_{a,x}^{(i)}|^2 1_{\{ |\xi_x^{(i)}|\leq a\}}\right)\\
 &=& \int_0^{a^2} \bbP(|\xi_x^{(i)}|>u^{\frac{1}{2}}){\rm d}u
\end{eqnarray*}
and this integral is equivalent as $a\rightarrow +\infty$ to $Ca^{2-\beta}$.\\
$\bullet$ Finally, the probability
$$\bbP\left(\hat \xi_{a,x}^{(i)} \neq 0 \right) = \bbP\left(|\xi_{x}^{(i)}| >a\neq 0 \right) $$
is equivalent to $(A_1+A_2)a^{-\beta}$ as $n\rightarrow +\infty$.

\noindent
{\it Proof of Proposition \ref{pr2}:}
For $a>0$, we use truncations with $a_n=a n^{\frac{1}{\alpha\beta}}c_n^{\frac{1}{\beta}}$ and write
$$G_n(t)=\bar \Gamma_{n,a}(t)+\hat\Gamma_{n,a}(t),$$
where
\begin{eqnarray*}
 \bar \Gamma_{n,a}(t)&=& n^{-\delta}c_n^{-\frac{1}{\beta}}\sum_{i=1}^{c_n}  \sum_{x\in\bbsZ} N_{nt}^{(i)}(x) \bar\xi_{a_n,x}^{(i)},\\
 \hat \Gamma_{n,a}(t)&=& n^{-\delta}c_n^{-\frac{1}{\beta}}\sum_{i=1}^{c_n}  \sum_{x\in\bbsZ} N_{nt}^{(i)}(x) \hat\xi_{a_n,x}^{(i)}.
\end{eqnarray*}
In order to prove Proposition \ref{pr2}, we prove that for any $T>0$,
\begin{equation}\label{eq3.19}
\lim_{a\rightarrow \infty}\limsup_{n\rightarrow\infty}\bbP\left( \sup_{t\in[0,T]}|\hat \Gamma_{n,a}(t)|>0\right)=0
\end{equation}
and that for fixed $a>0$, the family of processes $\left(\bar \Gamma_{n,a}(t)\right)_{t\geq 0}$ is tight.\\
Notice that the process $\hat G_n$ vanishes on $[0,T]$ as soon as the different remainder sceneries vanish on the range of the random walks up to time $[nT]+1$. Since the range of the $i-th$ random walk is included in $[-M_{[nT]+1}^{(i)},M_{[nT]+1}^{(i)}]$ and since the different random walks in random sceneries are independent, 
\begin{eqnarray*}
\bbP\left( \sup_{t\in[0,T]}|\hat \Gamma_{n,a}(t)|=0\right)&\geq& 
\left[ \bbP\left({\rm for\ all\ } x\in [-M_{[nT]+1}^{(1)},M_{[nT]+1}^{(1)}],\  \hat\xi_{a_n,x}^{(1)}=0 \right)\right]^{c_n}\\
 &\geq & \left(\bbE\left[\left(\bbP (\hat\xi_{a_n,0}^{(1)}=0)\right)^{2M_{[nT]+1}+1}\right]\right)^{c_n}
\end{eqnarray*}
Using Lemma \ref{lem4}, there exists a constant $C$ such that 
\begin{eqnarray*}
 \bbE\left[\left(\bbP (\hat\xi_{a_n,0}^{(1)}=0)\right)^{2M_{[nT]+1}+1}\right] &\geq& \bbE \left[\left(1-Ca_n^{-\beta} \right)^{2M_{[nT]+1}+1}\right]\\
 &\geq & \bbE \left[1+(2M_{[nT]+1}+1)\log(1-Ca_n^{-\beta}) \right]
\end{eqnarray*}
We majorate from above $\bbE(M_{[nT]+1})$ using equations (\ref{eq3.12}) and (\ref{eq1.02}), this yields
$$\bbE\left[\left(\bbP (\hat\xi_{a_n,0}^{(1)}=0)\right)^{2M_{[nT]+1}+1}\right]\geq 1+\log(1-Ca_n^{-\beta})C(nT)^{\frac{1}{\alpha}}.$$
Combinig these estimates, 
$$\bbP\left( \sup_{t\in[0,T]}|\hat \Gamma_{n,a}(t)|=0\right)\geq \left(1+\log(1-Ca_n^{-\beta})C(nT)^{\frac{1}{\alpha}}\right)^{c_n}. $$
Let $n$ go to infinity,
$$\limsup_{n\rightarrow\infty}\bbP\left( \sup_{t\in[0,T]}|\hat \Gamma_{n,a}(t)|>0\right)\leq 1-\exp(-Ca^{-\beta}T^{\frac{1}{\alpha}}).$$
Letting $a$ go to infinity, we obtain equation (\ref{eq3.19}).

For fixed $a$, we now prove the tightness of the family of processes $\left(\bar \Gamma_{n,a}(t)\right)_{t\geq 0}$.
Using the independence of the different random walks and random sceneries, we compute 
\begin{eqnarray}
& &\bbE\left[\left|\ \bar \Gamma_{n,a}(t_2)-\bar \Gamma_{n,a}(t_1)\ \right|^2\right]\nonumber\\
&=&{\rm Var}\left[\ \bar \Gamma_{n,a}(t_2)-\bar \Gamma_{n,a}(t_1)\  \right] +\left(\bbE\left[\ \bar \Gamma_{n,a}(t_2)-\bar \Gamma_{n,a}(t_1)\ \right]\right)^2\nonumber\\
&=&n^{-2\delta}c_n^{-\frac{2}{\beta}} \bbE\left[\sum_{i=1}^{c_n}\sum_{x\in\bbsZ} (N_{nt_2}^{(i)}(x)-N_{nt_1}^{(i)}(x))^2 \right]\bbE\left[\left(\bar\xi_{a_n,0}-\bbE(\bar\xi_{a_n,0}) \right)^2 \right]\nonumber\\
& &+n^{-2\delta}c_n^{-\frac{2}{\beta}} \left(\bbE\left[ \sum_{i=1}^{c_n}\sum_{x\in\bbsZ} (N_{nt_2}^{(i)}(x)-N_{nt_1}^{(i)}(x)) \right]\bbE\left[\bar\xi_{a_n,0}\right]\right)^2\label{eq3.20}
\end{eqnarray}
Using Lemma \ref{lem4},
there exists some $C$ such that 
$$|\bbE(\bar\xi_{a_n,0})|\leq Ca_n^{1-\beta}$$
and
$$\bbE\left[\left(\bar\xi_{a_n,0}-\bbE(\bar\xi_{a_n,0}) \right)^2 \right]\leq \bbE\left[\left(\bar\xi_{a_n,0} \right)^2 \right]\leq Ca_n^{2-\beta}.$$
Furthermore, 
$$\bbE\left[ \sum_{i=1}^{c_n}\sum_{x\in\bbsZ} (N_{nt_2}^{(i)}(x)-N_{nt_1}^{(i)}(x)) \right]=c_nn(t_2-t_1).$$
Using the strong Markov property and Lemma \ref{lem3.1}, the number of self-intersection between time $nt_1$ and $nt_2$ is overestimated by 
\begin{eqnarray*}
\bbE\left[\sum_{i=1}^{c_n}\sum_{x\in\bbsZ} (N_{nt_2}^{(i)}(x)-N_{nt_1}^{(i)}(x))^2 \right]&\leq&
\bbE\left[\sum_{i=1}^{c_n}\sum_{x\in\bbsZ} (N_{[nt_2]+1-[nt_1]}^{(i)}(x))^2 \right]\\
&\leq&  Cc_n([nt_2]-[nt_1]+1)^{2-\frac{1}{\alpha}}\\
\end{eqnarray*}
Furthermore, in the case  $|t_2-t_1|\leq 1/n$, we can see that  
$$\bbE\left[\sum_{i=1}^{c_n}\sum_{x\in\bbsZ} (N_{nt_2}^{(i)}(x)-N_{nt_1}^{(i)}(x))^2 \right]\leq 2c_n (nt_2-nt_1)^2,$$
since in the sum, at most $2c_n$ terms are not zero and those terms are bounded by $(nt_2-nt_1)^2$.\\
Combining equation all these estimates and equation (\ref{eq3.20}),
$$\bbE\left[\left|\ \bar \Gamma_{n,a}(t_2)-\bar \Gamma_{n,a}(t_1)\ \right|^2\right]\leq C|t_2-t_1|^2+C|t_2-t_1|^{2-\frac{1}{\alpha}}.$$
Using theorem 12.3 in Billingsley, this proves the tightness of the family of processes $\left(\bar \Gamma_{n,a}(t)\right)_{t\geq 0}$.
\\ $\Box$

\section{Properties of the process $\Gamma$ \\ Proof of Theorems \ref{theo3} \ref{theo4} \ref{theo5}}

In this section, we study the properties of the sample path of the limit process $\Gamma$.
The key of this study is a representation of the process $\Gamma$ as a stable integral (Theorem \ref{theo3}).
For the definition and properties of stable random measures and integrals with respect with these measures, we refer to the monography of Samorodnitsky and Taqqu \cite{ST94}. In \cite{CS05}, the authors study closely related processes, and we follow their methods. This allows us to prove various properties of the sample path of the limit process $\Gamma$ such as the probability tail of the suprema (Theorem \ref{theo4}) and the Hölder regularity (Theorem \ref{theo5}).
 
\subsection{Integral representation of $\Gamma$}
{\it Proof of Theorem \ref{theo3}:}\ \\
As we mentioned in a remark, equation (\ref{eq1.12}) is a consequence of the proof of Proposition \ref{pr1.1}. We now prove that $\Gamma$ has an integral representation given by equation (\ref{eq1.13}). \\
Let $G(t)$ be the process defined by the stable integral 
$$G(t)= \int_{\Omega'\times \bbsR} \sigma L_t(x)(\omega')M({\rm d}x,{\rm d}\omega').$$
where $M$ is a $\beta$-stable measure with control measure $\bbP'\times\lambda$ and skewness intensity $\nu$. To check that $G$ is well defined we need to prove that for all $t\geq 0$,
$$\int_{\Omega'\times \bbsR} \sigma^\beta L_t^{\beta}(x)(\omega'){\rm d}x\bbP'({\rm d}\omega')=\sigma^\beta \bbE'\left(\int_{\bbsR} L_t^{\beta}(x){\rm d}x \right) <\infty.$$
This was proven in Lemma \ref{lem2.1}, and hence the stable integral process $G$ is well defined. 

By definition of stable integration (see section 3 in \cite{ST94}), the finite dimensional distributions of the process $G$ are given by: \\
for any $\theta_1,\cdots,\theta_k \in\bbR$ and $t_1,\cdots,t_k\geq 0$
\begin{eqnarray*}
& &\bbE \exp\left(i\sum_{j=1}^k \theta_j G(t_j) \right)\\
&=& \exp\left(-\sigma^\beta\bbE'\int_{\bbsR}|\sum_{j=1}^k \theta_j L_{t_j}(x)|^\beta (1-i\nu\sgn(\sum_{j=1}^k \theta_j L_{t_j}(x)))\ {\rm d}x \right).
\end{eqnarray*}
We recognize the finite dimensional distribution of the process $\Gamma$ (see equation (\ref{eq1.12})). Hence, $G$ and $\Gamma$ have the same finite dimensional distributions and $G$ is an integral representation of the process $\Gamma$.
\\ $\Box$

\subsection{Probability tail of the maximum}
{\it Proof of Theorem \ref{theo4}: }\ \\
The integral representation of $\Gamma$ stated in Theorem \ref{theo3} is a key tool in the study of the path properties. A general theory is indeed available for the path properties of stable integral process.  The process $\Gamma$ has almost surely continuous path and hence is almost surely bounded on $[0,T]$, for every $T>0$.  We can hence apply Theorem 10.5.1 in \cite{ST94} giving the probability tail of suprema of bounded $\beta$-stable process with index $0<\beta<2$. This yields
\begin{eqnarray*}
& &\lim_{u\rightarrow \infty} u^{\beta}\bbP\left(\sup_{t\in [0,T]} \Gamma(t)\geq u \right)\\ 
&=&C_\beta\frac{1}{2} \left[ \bbE\left(\int_{\bbsR}L_T^+(x)^\beta(1+\nu) {\rm d}x\right)+\bbE\left(\int_{\bbsR}L_T^-(x)^\beta(1-\nu) {\rm d}x\right)\right],
\end{eqnarray*}
where $L_T^+$ and $L_T^-$ are defined by
$$L_T^+(x)=\sup_{0\leq t\leq T}\max(L_t(x),0)\ \ {\rm and}\ \ L_T^-(x)=\sup_{0\leq t\leq T}\max(-L_t(x),0).$$
We obtain the desired result since $L_T^+\equiv L_T$ and $L_T^-\equiv 0$.\\
The estimation of probability tail of $\sup_{t\in [0,T]} |\Gamma(t)|$ is proved in the same way. 
\\ $\Box$

\subsection{Hölder continuity of sample paths}

{\it Proof of Theorem \ref{theo5}: }\ \\
First consider the case $0<\beta<1$ and the process indexed by $(s,t)$ in $T=\{(s,t),0\leq s<t\leq 1/2\}$ defined by
$$\frac{G(t)-G(s)}{(t-s)^{1-\frac{1}{\alpha}}|\log(t-s)|^{\frac{1}{\alpha}}}=
\sigma^\beta\int_{\Omega'}\int_{\bbsR} \frac{ L_t(x)(\omega')- L_s(x)(\omega')}{(t-s)^{1-\frac{1}{\alpha}}|\log(t-s)|^{\frac{1}{\alpha}}}M({\rm d}x,{\rm d}\omega').$$
We have to prove that this process is sample bounded on $T$. Appliying Theorem 10.4.2 in \cite{ST94}, it is sufficient to check that 
\begin{equation}\label{eq3.21}
\sigma^\beta \int_{\Omega'}\int_{\bbsR}  \left(\sup_{(s,t)\in T}\frac{ |L_t^{\beta}(x)(\omega')- L_s^{\beta}(x)(\omega')|}{(t-s)^{1-\frac{1}{\alpha}}|\log(t-s)|^{\frac{1}{\alpha}}}\right)^\beta{\rm d}x\bbP'({\rm d}\omega')<\infty.
\end{equation}
Denote by 
$$K(x)=\sup_{(s,t)\in T}\frac{ |L_t^{\beta}(x)(\omega')- L_s^{\beta}(x)(\omega')|}{(t-s)^{1-\frac{1}{\alpha}}|\log(t-s)|^{\frac{1}{\alpha}}}.$$
and 
$$K=\sup_{x\in\bbsR}K(x).$$
It follows from the Hölder properties of local times of stable Lévy motions that $K$ is almost surly finite and has moments of all orders (see \cite{Ehm81}). Using Hölder inequality with $p>\frac{\alpha}{\alpha-1}$, the left hand side of eq (\ref{eq3.21}) is bounded from above
\begin{eqnarray*}
\bbE'\int_{\bbsR} K(x)^\beta {\rm d}x &\leq& \bbE'\int_{\bbsR} K^\beta 1_{\{\sup_{0\leq t\leq 1/2} |Y(t)|\geq |x|\}}{\rm d}x\\
&\leq & \bbE'(K^{\beta p})^{1/p}\int_{\bbsR}\bbP(\sup_{0\leq t\leq 1/2} |Y(t)|\geq |x|)^{1-1/p}{\rm d}x
\end{eqnarray*}
and this quantity is finite as soon as $(1-1/p)\alpha >1$ since $\bbE'(K^{\beta p})<\infty$ and $\bbP(\sup_{0\leq t\leq 1/2} |Y(t)|\geq |x|) \sim C|x|^{-\alpha}$ as $|x|\rightarrow\infty$.\\
It is worth nothing that this method yields an estimation on the tail probability :
$$\bbP\left( \sup_{0\leq s<t\leq 1/2}\frac{|\Gamma(t)-\Gamma(s)|}{(t-s)^{1-\frac{1}{\alpha}}|\log(t-s)|^{\frac{1}{\alpha}}}\geq \lambda\right)\sim C\lambda^{-\beta}$$
as $\lambda\rightarrow\infty$, with $C$ equal to the left hand side term of equaton (\ref{eq3.21}).

\noindent
In the case $1\leq \beta <2$, equation (\ref{eq3.21}) still holds, but this is not a sufficient condition for sample boundedness. Following Cohen and Samorodnitsky \cite{CS05}, we can use a series representation available in the symmetric case $\nu=0$, which shows that the process is conditionaly Gaussian and then use results on moduli of continuity of Gaussian process. The proof is the same, with local time of fractionnal Brownian motion replaced by local time of $\alpha$-stable Lévy motion.
$\Box$

\end{document}